\documentstyle[leqno]{article}

\newtheorem{th}{Theorem}[section]
\newtheorem{prop}[th]{Proposition}
\newtheorem{cor}[th]{Corollary}
\newtheorem{lema}[th]{Lemma}
\newcounter{defin}[section]
\renewcommand{\thedefin}{\thesection.\arabic{defin}}
\newcounter{ex}[section]

\newcounter{rem}[section]
\renewcommand{\therem}{\thesection.\arabic{rem}}
\newcommand{\T}{I\!\!\!T\!\!\!\!I}
\title{From PDE Systems and Metrics to\\
Generalized Field Theories}
\date{}
\author{Constantin UDRI\c STE and Mircea NEAGU}
\begin{document}
\maketitle
\begin{abstract}
Let $(T,h)$ and $(M,\varphi)$ be two Riemannian manifolds and
$$
(J^1(T,M),S=h+\varphi+h^{-1}*\varphi)
$$
the first-order jet fibre bundle,
endowed with Sasakian-like metric $S$, associated to these manifolds.
Developing our ideas from \cite{7}, \cite{22}, we show that a given
first-order {\bf PDEs} system on $J^1(T,M)$, and the Riemannian metric $S$,
determine an {\bf E}lectro{\bf D}ynamics {\bf M}etrical {\bf M}ulti-{\bf T}ime
{\bf L}agrange {\bf S}pace
$$
PDEsEDML^n_p=(J^1(T,M),L),
$$
where $L:J^1(T,M)\to R$ is a quadratic Lagrangian function of {\it
electrodynamics type} \cite{9}, In this new geometrical  structure, all $C^2$
solutions of the starting PDEs system become {\it harmonic maps},
being extremals of a least squares problem of variational calculus.
Our ideas are structured in the following way:

1) we find a suitable geometrical structure on $J^1(T,M)$ that convert the
solutions of a given PDEs system into harmonic maps (Section 1);

2) we build a natural geometry induced  by a  such PDEs system (Section 2);

3) we construct a field theory, in a general  setting, naturally attached
to this PDEs system (Section 3).

Consequently, we give a complete answer to our  generalization \cite{8},
\cite{20} of a problem rised by Poincar\'e, studied by others, but unfinalized
for a long time because of the absence of a suitable geometrical structure. The
Poincar\'e problem was solved recently  by the first author \cite{18}, using
the Lagrangian Geometry \cite{6}.

\end{abstract}

\section{Generalized Lorentz-Udri\c ste world-force law and a generalized
approach of inverse problem}

\hspace{5mm} Let  $T$  and $M$ be two smooth,  real, connected manifolds of
dimensions $p$ and $n$, with coordinates $(t^\alpha)_{\alpha=\overline{1,p}}$
and $(x^i)_{i=\overline{1,n}}$. Greek, respectively Latin, letters will be
used for  indexing the components of geometrical  objects attached to the
manifold $T$, respectively $M$.

From physical point of view, $T$ is  regarded  like a {\it "multi-temporal"}
manifold while the manifold $M$ like a  {\it "spatial"} one. Let us consider
the jet fibre bundle of order  one  $J^1(T,M)$, which is well known as a
basic object in the study  of classical  and quantum field theories.  The
local coordinates on $J^1(T,M)$ are $(t^\alpha,x^i,x^i_\alpha)$, where
$x^i_\alpha$ have the physical meaning of {\it "partial
directions"} or  {\it "partial derivatives"}.

Let us start with a given d-tensor  $X^{(i)}_{(\alpha)}(t^\gamma,x^k)$ on
$J^1(T,M)$, which defines the PDEs system of order one
\begin{equation}\label{pde}
x^i_\alpha=X^{(i)}_{(\alpha)}(t^\gamma,x^k(t^\gamma)),
\end{equation}
where $\displaystyle{x^i_\alpha={\partial x^i\over\partial t^\alpha}}$.
Obviously, the complete integrability conditions,
$$
{\partial  X^{(i)}_{(\alpha)}\over\partial t^\beta}+{\partial  X^{(i)}_{(
\alpha)}\over\partial x^m}X^{(m)}_{(\beta)}={\partial  X^{(i)}_{(\beta)}\over
\partial t^\alpha}+{\partial  X^{(i)}_{(\beta)}\over\partial x^m}X^{(m)}_{(
\alpha)},
$$
are required by the existence of solutions.

The Lorentz world-force law, initially  stated   for particles
in non-quantum relativity \cite{14}, was generalized by Udri\c ste, using
the notion of  {\it potential map} \cite{20}, \cite{21}. In this direction,
starting with $h=(h_{\alpha\beta}(t^\gamma))$  and $\varphi=(\varphi_{ij}(x^k))$
as semi-Riemannian metrics on $T$, respectively $M$, the first author proved
the following

\begin{th}{(Lorentz-Udri\c ste World-Force Law)}\label{lu}\medskip\\
Every  $C^2$ solution of the PDEs system $(1.1)$ is a horizontal potential
map of the  semi-Riemannian-Lagrange  manifold
$$
(T\times  M,\; h+\varphi,\; M^{(i)}_{(\alpha)\beta}=-H^\mu_{\alpha\beta}x^i_\mu,
\;N^{(i)}_{(\alpha)j}=\gamma^i_{jm}x^m_\alpha-F^{\;\;i}_{j\alpha}),
$$
where $H^\mu_{\alpha\beta}$ (resp. $\gamma^i_{jk}$) are the Christoffel symbols
of  $h_{\alpha\beta}$ (resp. $\varphi_{ij}$), and
$$
F^{\;\;i}_{j\alpha}=\nabla_jX^{(i)}_{(\alpha)}-\varphi^{ir}\varphi_{js}\nabla_rX^{(s)}
_{(\alpha)}
$$
is the external distinguished tensor  field which characterizes the
helicity of the distinguished tensor field $X^{(i)}_{(\alpha)}$.
\end{th}

A new version of  the Lorentz-Udri\c ste World-Force  Law was formulated
by the second author. Its advantage is the possibility of construction a
generalized field theory, naturally attached to the PDEs
system $(1.1)$.

In order to give  this version, we suppose that the metrics $h$ and $\varphi$
are just Riemannian metrics and we endow the jet fibre bundle $J^1(T,M)$
with the  canonical nonlinear connection $\Gamma_0=(\stackrel{0}{M}\;\!\!\!\!^{(i)}_
{(\alpha)\beta},\stackrel{0}{N}\;\!\!\!\!^{(i)}_{(\alpha)j})$ associated to the
Riemannian metrics pair $(h,\varphi)$, where
$$
\begin{array}{l}\medskip
\stackrel{0}{M}\;\!\!\!^{(i)}_{(\alpha)\beta}=-H^\mu_{\alpha\beta}x^i_\mu,\\
\stackrel{0}{N}\;\!\!\!^{(i)}_{(\alpha)j}=\gamma^i_{jm}x^m_\alpha.
\end{array}
$$
At the same  time,  let us consider $B\Gamma_0=(H^\gamma_{\alpha\beta},0,
\gamma^i_{jk},0)$ the Berwald  $\Gamma_0$-linear connection on $J^1(T,M)$ and
$"_{//\beta}"$, $"_{\Vert j}"$,  $"\Vert^{(\beta)}_{(i)}"$ its local covariant
derivatives \cite{10}.

Taking into account the $x^m_\mu$-independence of the
d-tensor field $X^{(i)}_{(\alpha)}$, by a simple calculation, we obtain
\begin{prop}
The horizontal local covariant derivatives of $X^{(i)}_{(\alpha)}$ have the
following expressions:
$$
\begin{array}{l}\medskip
\displaystyle{X^{(i)}_{(\alpha)//\beta}={\partial X^{(i)}_{(\alpha)}\over
\partial t^\beta}-X^{(i)}_{(\mu)}H^\mu_{\alpha\beta}},\\
\displaystyle{X^{(i)}_{(\alpha)\Vert j}={\partial X^{(i)}_{(\alpha)}\over
\partial x^j}+X^{(m)}_{(\alpha)}\gamma^i_{mj}}.
\end{array}
$$
\end{prop}

In this context, supposing that $T$ is a compact and orientable manifold,
it follows the following generalization of the Theorem \ref{lu}.

\begin{th}\label{lun}(Generalized Lorentz-Udri\c ste World-Force Law)\medskip\\
The $C^2$ solutions of the PDEs system $(1.1)$ are harmonic maps
of the multi-time dependent spray $(H,G)$ on $J^1(T,M)$, defined by
the temporal components
\begin{equation}\label{ts}
H^{(i)}_{(\alpha)\beta}=-{1\over 2}H^\gamma_{\alpha\beta}x^i_\gamma,
\end{equation}
and the local spatial components
$$
G^{(i)}_{(\alpha)\beta}={1\over 2}\gamma^i_{jk}x^j_\alpha x^k_\beta+h_{\alpha
\beta}F^i,
$$
where,
$$
F^i={h^{\mu\nu}\over 2p}\left\{\varphi^{il}X^{(s)}_{(\nu)\Vert l}\varphi_{sr}
\left[X^{(r)}_{(\mu)}-x^r_\mu\right]+X^{(i)}_{(\nu)\Vert m}x^m_\mu+X^{(i)}
_{(\nu)//\mu}\right\},\quad p=\dim T.
$$
In other words, the $C^2$ solutions of $(1.1)$ verify the harmonic map
equations \cite{9}, \cite{11},
$$
h^{\alpha\beta}\left\{x^i_{\alpha\beta}+2H^{(i)}_{(\alpha)\beta}+2G^{(i)}_
{(\alpha)\beta}\right\}=0.
$$
\end{th}
{\bf Proof.} Let  us  consider the multi-time least squares Lagrangian
$$
\begin{array}{l}\medskip
{\cal L}=\Vert\hbox{\bf C}-\hbox{\bf X}\Vert^2\sqrt{h}=\left\{h_{\alpha\beta}
(t^\gamma)\varphi_{ij}(x^k)\left[x^i_\alpha-X^{(i)}_{(\alpha)}\right]
\left[x^j_\beta-X^{(j)}_{(\beta)}\right]\right\}\sqrt{h}=\\
\hspace*{4mm}=\left\{h^{\alpha\beta}(t^\gamma)\varphi_{ij}(x^k)x^i_\alpha
x^j_\beta+U^{(\alpha)}_{(i)}(t^\gamma,x^k)x^i_\alpha+\Phi(t^\gamma,x^k)\right\}
\sqrt{h},
\end{array}
$$
where
\begin{equation}\label{rel}
\begin{array}{l}\medskip
\displaystyle{\hbox{\bf C}=x^i_\alpha{\partial\over\partial x^i_\alpha},}
\quad
U^{(\alpha)}_{(i)}=-2h^{\alpha\mu}\varphi_{im}X^{(m)}_{(\mu)},
\\
\displaystyle{\hbox{\bf X}=X^{(i)}_{(\alpha)}{\partial\over\partial x^i_\alpha},}
\quad
\Phi=h^{\mu\nu}\varphi_{rs}X^{(r)}_{(\mu)}X^{(s)}_{(\nu)},
\end{array}
\end{equation}
and let ${\cal E}_{\cal L}:C^2(T,M)\to R_+$ (note that the Riemannian
metrics $h_{\alpha\beta}$ and $\varphi_{ij}$ are  positive definite) be the
least squares energy of ${\cal L}$, i. e. ,
$$
{\cal E}_{\cal L}=\int_T{\cal L}dt^1\wedge dt^2\wedge\ldots\wedge dt^p=
\int_T\Vert\hbox{\bf C}-\hbox{\bf X}\Vert^2\sqrt{h}dt^1\wedge dt^2\wedge
\ldots\wedge dt^p\geq 0.
$$

It is obvious now that $f\in C^2(T,M)$, locally expressed by $(t^\gamma)\to
(x^i(t^\gamma))$, is a solution of the PDEs system $(1.1)$ iff the  map  $f$
is a global minimum point for the least squares  energy functional ${\cal E}_{\cal L}$.
Therefore, a $C^2$ solution $f$ of $(1.1)$ verifies the Euler-Lagrange
equations,
\begin{equation}\label{el}
{\partial{\cal L}\over\partial x^i}-{\partial\over\partial t^\alpha}\left(
{\partial{\cal L}\over\partial x^i_\alpha}\right)=0,\quad\forall\;i=\overline
{1,n}.
\end{equation}
But, we proved in the paper  \cite{9} that the Euler-Lagrange
equations $(1.4)$ can be rearranged in the general Poisson form,
$$
h^{\alpha\beta}\left\{x^i_{\alpha\beta}+2H^{(i)}_{(\alpha)\beta}+2G^{(i)}_
{(\alpha)\beta}\right\}=0,
$$
where the temporal spray components $H^{(i)}_{(\alpha)\beta}$ are  given by
$(1.2)$, while the spatial spray components have the expressions
$$
G^{(i)}_{(\alpha)\beta}={1\over 2}\gamma^i_{jk}x^j_\alpha x^k_\beta+
{h_{\alpha\beta}g^{il}\over 4p}\left[U^{(\mu)}_{(l)m}x^m_\mu+{\partial U^{(\mu)}
_{(l)}\over\partial t^\mu}+U^{(\mu)}_{(l)}H^\gamma_{\mu\gamma}-{\partial \Phi
\over\partial x^l}\right],
$$
where $p=\dim T$ and
$$
\displaystyle{U^{(\alpha)}_{(i)j}={\partial U^{(\alpha)}_{(i)}\over\partial
x^j}-{\partial U^{(\alpha)}_{(j)}\over\partial x^i}}.
$$
Now, using  the relations $(1.3)$ and direct computations, it follows
$$
\begin{array}{l}\medskip
U^{(\alpha)}_{(i)j}=-2h^{\alpha\mu}\left[\varphi_{im}X^{(m)}_{(\mu)\Vert j}-
\varphi_{jm}X^{(m)}_{(\mu)\Vert i}\right],\\\medskip
\displaystyle{{\partial U^{(\alpha)}_{(i)}\over\partial t^\alpha}+U^{(\nu)}_
{(i)}H^\alpha_{\nu\alpha}=-2h^{\mu\nu}\varphi_{im}X^{(m)}_{(\mu)//\nu}},\\
\displaystyle{{\partial\Phi\over\partial x^l}=2h^{\mu\nu}\varphi_{mr}X^{(m)}_
{(\mu)}X^{(r)}_{(\nu)\Vert l}}.\;\;\mbox{\rule{5pt}{5pt}}
\end{array}
$$
\addtocounter{rem}{1}
{\bf Remark \therem} The method used in the proof of the theorem  \ref{lun}
is a general one and may be  called the {\it least squares variational  calculus}
method. For   example, let us consider $X$ a d-tensor field on the jet fibre bundle
of $r$-order $J^r(T,M)$, supposed as being endowed with an {\it "a priori"}
Riemannian metric $<,>$, and  let $X(t,x^{(r)})=0$ be the associated $r$-order
PDEs system. Automatically, we can attach a least squares problem of variational
calculus, using the Lagrangian function $L=\Vert X\Vert^2$, where the norm
$\Vert\;\Vert$ is buildet by the metric $<,>$.  This remark works  on any
DEs or PDEs  system. As example, the Yang-Mills functional
$$
{\cal YM}(\nabla)={1\over 2}\int_M\Vert R^\nabla\Vert^2d{\cal V}_g,
$$
on the space of connections $\nabla$, is the least squares functional attached
to the PDEs system $R^\nabla=0$. In this sense, ${\cal  YM}(\nabla)$ measures the
deviation from flatness.\medskip

The theorem \ref{lun}  shows that the $C^2$ solutions of the PDEs system
$(1.1)$ can be viewed like some {\it harmonic maps} on $J^1(T,M)$,
via two {\it "a  priori"} fixed Riemannian metrics $h$ on $T$, respectively
$\varphi$ on $M$.

We recall the well known fact that paracompactness of a real manifold
implies the existence of a Riemannian metric. So, we can restrict the imposed
condition concerning the {\it "a priori"} existence of the Riemannian metrics
$h$ and $\varphi$, supposing the paracompactness of the manifolds $T$
and $M$.

In conclusion, we can formulate the following
\begin{cor}
If $T$ and $M$ are paracompact manifolds, then the $C^2$ solutions of the
PDEs system of order one $(1.1)$ can be regarded  like harmonic maps on the
jet fibre  bundle of order one $J^1(T,M)$, in the sense of theorem \ref{lun}.
\end{cor}

In the sequel, let us denote $J^s(T,M)$ the  total  space of the jet fibre
bundle of order $s$, whose local coordinates are
$$
(t,x^{(s)})=(t^\gamma,x^k,x^k_{\gamma_1},x^k_{{\gamma_1}{\gamma_2}},\ldots,x^k_{{
\gamma_1}{\gamma_2}\ldots{\gamma_s}}),
$$
where the  coordinates $x^k_{{\gamma_1}{\gamma_2}\ldots{\gamma_l}},\;l\in\{
1,2,\ldots ,s\}$, $k\in\{1,2,\ldots ,n\}$,  have the meaning of {\it "partial
derivatives of order $l$ of the spatial variables $x^k$, with respect to the
temporal variables $t^{\gamma_1},t^{\gamma_2},\ldots,t^{\gamma_l}$".}
\begin{lema}
The dimension of the total space of the jet fibre bundle $J^s(T,M)$ is
$$
p+n\left[C^0_p+C^1_p+\ldots+C^{s-1}_p\right],
$$
where $\dim T=p$ and $\dim M=n$.
\end{lema}

Let us consider the PDEs system of order $r\geq 1$ on $J^r(T,M)$, locally
expressed by
\begin{equation}\label{pder}\hspace*{5mm}
x^i_{\alpha_1\alpha_2\ldots\alpha_{r-1}\alpha_r}=X^{(i)}_{\alpha_1\alpha_2
\ldots\alpha_{r-1}(\alpha_r)}(t,x^{(r-1)}),
\end{equation}
where $X^{(i)}_{\alpha_1\alpha_2\ldots\alpha_{r-1}(\alpha_r)}$ is a
d-tensor field on $J^r(T,M)$, with respect to the indices $i$ and $\alpha_r$.

Let $\tilde M$ be  the  submanifold of $J^r(T,M)$,  whose coordinates are
only \linebreak
$x^{(r-1)}=(x^k,x^k_{\gamma_1},x^k_{{\gamma_1}{\gamma_2}},\ldots,x^k_{{\gamma_1}
{\gamma_2}\ldots{\gamma_{r-1}}})$. In these conditions, we deduce the
following  important result.

\begin{th}{(Geometrical Solution of Generalized Inverse Problem)}\medskip\\
If the manifolds $T$  and $M$ are paracompact, then the $C^{r+1}$ solutions
of the PDEs system $(1.5)$ can be viewed like harmonic maps on the
jet fibre bundle of order one $J^1(T,\tilde M)$, whose dimension is
$$
p+(p+1)\left[C^0_p+C^1_p+\ldots+C^{r-1}_p\right]n,
$$
where $\dim T=p$ and $\dim M=n$.
\end{th}
{\bf Proof.} Obviously, the paracompactness of the manifolds $T$ and $M$
implies  the  paracompactness of the  jet fibre bundle $J^r(T,M)$ and,
implicitly the paracompactness  of  $\tilde M$. Moreover, if
$(x^{(r-1)})=(x^k,x^k_{\gamma_1},x^k_{{\gamma_1}{\gamma_2}},\ldots,x^k_{{\gamma_1}
{\gamma_2}\ldots{\gamma_{r-1}}})$ are the  coordinates of $\tilde M$ and
$(x^{i(r-1)})=(x^i_{{\alpha_1}{\alpha_2}\ldots {\alpha_{r-1}}})$, then the
PDEs system $(1.5)$  writes in the form,
$$
x^{i(r-1)}_{\alpha_r}=X^{(i)}_{(\alpha_r)}(t,x^{(r-1)}).
$$
\nopagebreak
Now, the theorem \ref{lun} implies the   required result. \rule{5pt}{5pt}
\pagebreak\\
\addtocounter{rem}{1}
{\bf Remark \therem} In a general approach of the well known inverse problem,
we consider that the previous theorem gives a  final geometrical answer for this
open problem. From our point of  view, the old classical open problem was not
solved because  of two reasons:

i) Firstly, it was incorrectly formulated;

ii)  Secondly, almost all attempts of solving this problem used an unsuitable
Lagrangian function space of possible solutions.

In our opinion, the inverse problem must be formulated in the following way:
\medskip\\
{\bf\underline{Generalized Inverse Problem:}} Find a space of Lagrangian functions such
that the $C^2$ solutions of a givenPDEs system to  become  solutions of the
Euler-Lagrange equations of a Lagrangian on this space.

\section{From PDEs system of order one and metrics to generalized Maxwell
and Einstein equations}

\setcounter{equation}{0}
\hspace{5mm} Let us return to the PDEs system of order one, given by $(1.1)$.
The theorem \ref{lun} ensures us that a $C^2$ map  $f$ is a solution of PDEs
system $(1.1)$ iff  $f$ is a global minimum point of the energy functional
produced by the multi-time Lagrangian ${\cal L}=L\sqrt{h}$, whose least
squares Lagrangian
function
$L:J^1(T,M)\to R$ is  expressed by
\begin{equation}\label{led}
L=\Vert\hbox{\bf C}-\hbox{\bf X}\Vert^2=h^{\alpha\beta}(t^\gamma)\varphi_{ij}(x^k)
x^i_\alpha x^j_\beta+U^{(\alpha)}_{(i)}(t^\gamma,x^k)x^i_\alpha+\Phi(t^\gamma,x^k),
\end{equation}
where  $U^{(\alpha)}_{(i)}=-2h^{\alpha\mu}\varphi_{im}X^{(m)}_{(\mu)}$ and
$\Phi=h^{\mu\nu}\varphi_{rs}X^{(r)}_{(\mu)}X^{(s)}_{(\nu)}$.

But, the  Lagrangian  function $L$ is a natural generalization  of the classical
Lagrangian of  electrodynamics \cite{6},  which governs the movement
law of a particule placed concomitantly into a gravitational field and an
electromagnetic one. For that reason, in our generalized  context, the Riemannian
metric $h_{\alpha\beta}$  (resp. $\varphi_{ij}$)  will have the {\it temporal}
(resp. {\it spatial}) {\it gravitational potential} physical meaning, the
components $X^{(\alpha)}_{(i)}$ that of {\it electromagnetic potentials} on
$J^1(T,M)$, and $F$  that of the {\it potential function}.

In this geometrical-physical context, following the terminology used in
\cite{9}, we can introduce
\medskip\\
\addtocounter{defin}{1}
{\bf Definition \thedefin}
The pair $(J^1(T,M),L)$, which consists of
jet fibre bundle of order one and a Lagrangian function of the form $(2.1)$,
is called the {\it canonical autonomous metrical multi-time Lagrange space
of electrodynamics, attached to the PDEs system $(1.1)$ and to the Riemannian
metrics $h$ and $\varphi$}, and is denoted $PDEsEDML^n_p$.\medskip

On the jet  fibre  bundle of order one $J^1(T,M)$, a natural generalized
theory of
physical fields, attached to a given multi-time Lagrangian  function of
electrodynamics type was recently developed by the second author \cite{9}.
Consequently, via the canonical autonomous  metrical multi-time Lagrange
space of electrodynamics $PDEsEDML^n_p$, we can construct a natural physical
field theory, in a general setting, arising only from the PDEs system $(1.1)$
and two {\it "a priori"}  given Riemannian metrics $h_{\alpha\beta}$ and
$\varphi_{ij}$.\medskip\\
\addtocounter{defin}{1}
{\bf Definition \thedefin} The multi-time dependent spray $(H,G)$ used
in the theorem \ref{lun} is called the {\it canonical multi-time  dependent
spray attached to $PDEsEDML^n_p$}.
\medskip

Using the canonical multi-time dependent spray $(H,G)$, one naturally induces
a nonlinear connection
$\Gamma=(M^{(i)}_{(\alpha)\beta}, N^{(i)}_{(\alpha)j})$ on
$J^1(T,M)$, which is also called the {\it  canonical nonlinear connection
of the autonomous metrical multi-time Lagrange space of electrodynamics
$PDEsEDML^n_p$.}
In this direction, denoting ${\cal G}^i=h^{\alpha\beta}G^{(i)}_{(\alpha)\beta}$,
we establish the folowing
\begin{th}
The canonical nonlinear connection of the autonomous metrical multi-time Lagrange
space of electrodynamics $PDEsEDML^n_p$ is determined by the temporal
components
$$
M^{(i)}_{(\alpha)\beta}=2H^{(i)}_{(\alpha)\beta}=\stackrel{0}{M}\;\!\!\!^{(i)}_
{(\alpha)\beta}
$$
and the local spatial components
$$
N^{(i)}_{(\alpha)j}={\partial{\cal G}^i\over\partial x^j_\gamma}h_{\alpha\gamma}
=\stackrel{0}{N}\;\!\!\!^{(i)}_{(\alpha)j}-F^{\;\;i}_{j\alpha},
$$
where $\stackrel{0}{M}^{(i)}_{(\alpha)\beta}=-H^\gamma_{\alpha\beta}x^i_
\gamma$, $\;\stackrel{0}{N}^{(i)}_{(\alpha)j}=\gamma^i_{jk}x^k_\alpha$, and
the d-tensor field
$$
F^{\;\;i}_{j\alpha}={1\over 2}\left[X^{(i)}_{(\alpha)\Vert j}-\varphi^{ir}
X^{(s)}_{(\alpha)\Vert r}\varphi_{sj}\right]
$$
characterizes the helicity  of the d-tensor field $X^{(i)}_{(\alpha)}$.
\end{th}

To determine the generalized {\it Cartan canonical connection} $C\Gamma$ of
the autonomous metrical multi-time Lagrange space of electrodynamics
$PDEsEDML^n_p$, together with its torsion and curvature local d-tensors,
we use the adapted dual  bases
$\displaystyle{\left\{{\delta\over\delta t^\alpha},
{\delta\over\delta x^i}, {\partial\over\partial x^i_\alpha}\right\}
\subset{\cal X}(E)}$ and $\{dt^\alpha, dx^i, \delta x^i_\alpha\}\subset
{\cal X}^*(E)$ of the nonlinear connection $\Gamma$, via the formulas
\begin{equation}
\begin{array}{l}\medskip
\displaystyle{{\delta\over\delta t^\alpha}={\partial\over\partial t^\alpha}-
M^{(j)}_{(\beta)\alpha}{\partial\over\partial x^j_\beta}}\\\medskip
\displaystyle{{\delta\over\delta x^i}={\partial\over\partial x^i}-
N^{(j)}_{(\beta)i}{\partial\over\partial x^j_\beta}}\\
\delta x^i_\alpha=dx^i_\alpha+M^{(i)}_{(\alpha)\beta}dt^\beta+N^{(i)}_{(\alpha)
j}dx^j.
\end{array}
\end{equation}

In this context, by simple  computations, we find
\begin{th}
i) The generalized Cartan canonical connection $C\Gamma$
of the metrical multi-time Lagrange space $PDEsEDML^n_p$
has the adapted components
$$
H^\gamma_{\alpha\beta}=H^\gamma_{\alpha\beta},\quad G^k_{j\gamma}=0, \quad
L^i_{jk}=\gamma^i_{jk},\quad C^{i(\gamma)}_{j(k)}=0.
$$

ii) The torsion {\em\bf T} of the generalized Cartan canonical connection
of the autonomous metrical multi-time Lagrange space of electrodynamics
$PDEsEDML^n_p$ is
determined by three adapted local d-tensors, namely,
$$
\begin{array}{l}\medskip
\displaystyle{
R^{(i)}_{(\alpha)\beta\gamma}=-H^\mu_{\alpha\beta\gamma}x^i_\mu,}
\\\medskip
\displaystyle{R^{(i)}_{(\alpha)\beta j}={1\over 2}\left[X^{(i)}_{(\alpha)
\Vert j//\beta}-\varphi^{ir}X^{(s)}_{(\alpha)\Vert r//\beta}\varphi_{sj}
\right],}\\
\displaystyle{
R^{(i)}_{(\alpha)jk}=r^i_{jkm}x^m_\alpha-{1\over 2}\left[X^{(i)}_{(\alpha)\Vert j
\Vert k}-\varphi^{ir}X^{(s)}_{(\alpha)\Vert r\Vert k}\varphi_{sj}\right],}
\end{array}
$$
where $H^\mu_{\alpha\beta\gamma}$ (resp. $r^l_{ijk}$) are the local curvature
tensors of the Riemannian metric $h_{\alpha\beta}$ (resp. $\varphi_{ij}$),
and
$$
\begin{array}{l}\medskip
\displaystyle{X^{(i)}_{(\alpha)\Vert j//\beta}={\partial X^{(i)}_{(\alpha)\Vert j}
\over\partial t^\beta}-X^{(i)}_{(\mu)\Vert j}H^\mu_{\alpha\beta}},\\
\displaystyle{X^{(i)}_{(\alpha)\Vert j\Vert k}={\partial X^{(i)}_{(\alpha)\Vert j}
\over\partial x^k}+X^{(m)}_{(\alpha)\Vert j}\gamma^i_{mk}-X^{(i)}_{(\alpha)
\Vert m}\gamma^m_{jk}}.
\end{array}
$$

iii) The curvature {\em\bf R} of the generalized Cartan canonical connection
of the autonomous
metrical multi-time Lagrange space of electrodynamics $PDEsEDML^n_p$ is
determined by two adapted local d-tensors, namely, $H^\delta_{\alpha\beta\gamma}$ and
$R^l_{ijk}=r^l_{ijk}$, that is, exactly the curvature tensors of the Riemannian
metrics $h_{\alpha\beta}$ and $\varphi_{ij}$.
\end{th}
\addtocounter{rem}{1}
{\bf Remarks \therem} i) Both Cartan  and Berwald  connections on
$J^1(T,M)$ have the same adapted components. However, they are two distinct
linear connections. This is because  the Cartan connection  is a $\Gamma$-linear
connection while the Berwald connection is a $\Gamma_0$-linear one.

ii) Let $"_{/\beta}"$, $"_{\vert j}"$ and $"\vert^{(\beta)}_{(j)}"$ be the
local covariant derivatives of the Cartan connection $C\Gamma$, and let us
consider a distinguished tensor field $D=(D^{\alpha i(j)(\delta)\ldots}_
{\gamma k(\beta)(l)\ldots})$, whose  components do not depend on the partial
directions $x^m_\mu$. In these conditions, we obtain without difficulties
that
$$
D^{\alpha i(j)(\delta)\ldots}_{\gamma k(\beta)(l)\ldots!A}=D^{\alpha i(j)
(\delta)\ldots}_{\gamma k(\beta)(l)\ldots!!A},
$$
where $"_{!A}"$, respectively $"_{!!A}"$, is one of the local operators
$"_{/\beta}"$, $"_{\vert j}"$ or $"\vert^{(\beta)}_{(j)}"$, respectively
$"_{//\beta}"$, $"_{\Vert j}"$ or $"\Vert^{(\beta)}_{(j)}"$.

iii)  The only geometrical objects of $PDEsEDML^n_p$, effectively
dependent of the PDEs system $(1.1)$, are the {\it spatial  nonlinear
connection components} $N^{(i)}_{(\alpha)j}$ and the {\it adapted local
torsion d-tensors} $R^{(i)}_{(\alpha)\beta j}$, $R^{(i)}_{(\alpha)jk}$.
\medskip

In the sequel, following  the paper \cite{9}, we shall write the generalized
Maxwell and Einstein equations produced by $PDEsEDML^n_p$.
\begin{th}\label{me}
i) The electromagnetic distinguished  2-form of $PDEsEDML^n_p$ is
$$
F=F^{(\alpha)}_{(i)j}\delta x^i_\alpha\wedge dx^j,
$$
where
$$
F^{(\alpha)}_{(i)j}={h^{\alpha\mu}\over 2}\left[\varphi_{im}X^{(m)}_{(\mu)
\Vert j}-\varphi_{jm}X^{(m)}_{(\mu)\Vert i}\right].
$$

ii) The electromagnetic local components $F^{(\alpha)}_{(i)j}$ of the autonomous
metrical multi-time Lagrange space of electrodynamics $PDEsEDML^n_p$ are
governed by the following generalized  Maxwell equations,
$$
\left\{\begin{array}{l}\medskip
\displaystyle{
F^{(\alpha)}_{(i)j//\beta}={1\over 4}{\cal A}_{\{i,j\}}\left\{h^{\alpha\mu}
\varphi_{im}\left[X^{(m)}_{(\mu)\Vert j//\beta}-\varphi^{mr}X^{(s)}_{(\mu)
\Vert r//\beta}\varphi_{sj}\right]\right\}}\\\medskip
\sum_{\{i,j,k\}}F^{(\alpha)}_{(i)j\Vert k}=0\\
\sum_{\{i,j,k\}}F^{(\alpha)}_{(i)j}\Vert^{(\gamma)}_{(k)}=0,
\end{array}\right.
$$
where ${\cal A}_{\{i,j\}}$ represents an alternate sum and $\sum_{\{i,j,k\}}$
means a cyclic sum.
\end{th}

\begin{th}\label{ee}
i) The gravitational $h$-potential on $J^1(T,M)$, induced by the metrical
multi-time Lagrange space $PDEsEDML^n_p$, is given by the Sasakian-like
metric
$$
G=h_{\alpha\beta}dt^\alpha\otimes dt^\beta+\varphi_{ij}dx^i\otimes dx^j+
h^{\alpha\beta}\varphi_{ij}\delta x^i_\alpha\otimes\delta x^j_\beta.
$$
ii) The generalized  Einstein equations which govern the gravitational
$h$-potential $G$ of $PDEsEDML^n_p$ are locally expressed by
$$
\left\{\begin{array}{l}\medskip
\displaystyle{H_{\alpha\beta}-{H+r\over 2}h_{\alpha\beta}={\cal K}{\cal T}_
{\alpha\beta}}\\\medskip
\displaystyle{r_{ij}-{H+r\over 2}\varphi_{ij}={\cal K}{\cal T}_{ij}}\\\medskip
\displaystyle{-{H+r\over 2}h^{\alpha\beta}\varphi_{ij}={\cal K}{\cal T}^{(\alpha)
(\beta)}_{(i)(j)}},
\end{array}\right.\leqno{(E_1)}
$$
$$
\left\{\begin{array}{lll}\medskip
0={\cal T}_{\alpha i},&0={\cal T}_{i\alpha},&
0={\cal T}^{(\alpha)}_{(i)\beta}\\
0={\cal T}^{\;(\beta)}_{\alpha(i)},&0={\cal T}^{\;(\alpha)}_{i(j)},&
0={\cal T}^{(\alpha)}_{(i)j},
\end{array}\right.\leqno{(E_2)}
$$
where ${\cal T}_{AB},\; A,B\in\{\alpha,i,{(\alpha)\atop(i)}\}$ represent the
components of an intrinsic stress-energy d-tensor ${\cal T}$ of matter
on $J^1(T,M)$, $H_{\alpha\beta}$ (resp.  $r_{ij}$)  are the components of
Ricci tensor  of $h_{\alpha\beta}$ (resp. $\varphi_{ij}$), and $H$ (resp. $r$)
is the scalar curvature  of the  Riemannian metric  $h_{\alpha\beta}$ (resp.
$\varphi_{ij}$).

iii) The components ${\cal T}_{\mu\beta}$ and  ${\cal T}_{mj}$ of the
stress-energy d-tensor ${\cal T}$ of the metrical multi-time Lagrange space
$PDEsEDML^n_p$ satisfy the classical conservation laws,
$$
{\cal T}^\mu_{\beta//\mu}=0,\quad
{\cal T}^m_{j\Vert m}=0.
$$
\end{th}
\addtocounter{rem}{1}
{\bf Remark \therem} The theorems \ref{me} and \ref{ee} emphasize that  only
the {\it electromagnetic field} $F$ of the $PDEsEDML^n_p$ is {\it  effectively}
dependent by the PDEs system $(1.1)$,
via  the {\it electromagnetic  potentials} $X^{(\alpha)}_{(i)}$. Consequently, in
order to  obtain some  information  upon the PDEs system $(1.1)$, via the
canonical attached generalized field theory, it is interesting to  study only
the {\it electromagnetism} generated by this system.

\section{Generalized electromagnetic theory induced by remarkable
PDE systems}

\setcounter{equation}{0}
\hspace{5mm} In this section we will consider some important particular PDE
systems of order  one and we apply them the previous general geometrical and
physical results.  Taking into account the remarks 2.1 and  2.2, we will study
only the geometrical and physical objects which are {\it effectively} determined
by these systems (i.e. , the nonlinear  connection, the torsion d-tensors and
the electromagnetic field). From this point of view, we will ignore their
attached generalized  gravitational theory, because it is dependent only
of the Riemannian metrics used in the construction of the metrical multi-time
Lagrange space $PDEsEDML^n_p$.

\subsection{Orbits}

\hspace{5mm} Let $T=[a,b]\subset  R,\;h_{11}(t)=1$ and $X^{(i)}_{(1)}(t,
x^k)=\xi^i(x^k)$, where $\xi^i(x^k)$ is a d-vector field  on $J^1([a,b],M)$.
The PDEs system $(1.1)$ becomes
$$
{dx^i\over dt}=\xi^i(x^k(t)),
\leqno{(O)}
$$
that is, the DEs system which gives the orbits  of the spatial vector field
$\xi$.

In this context, the $C^2$ orbits of $\xi$  are the global minimum  points of
the least squares energy  functional
$$
{\cal E}_\xi=\int^b_a\left[\varphi_{ij}(x^k)y^iy^j-2\xi_j(x^k)y^j+\Vert\xi\Vert^2
(x^k)\right]dt,
$$
where $y^i=x^i_1=dx^i/dt,\;\xi_j(x^k)=\varphi_{jm}(x^k)\xi^m(x^k)$ and $\Vert\xi\Vert^2
(x^k)=\xi^m(x^k)\xi_m(x^k)$.

\begin{th}
i) The canonical nonlinear  connection, induced by the
system $(O)$ on $J^1([a,b],M)\equiv TM$, is given by the components
$$M^{(i)}_{(1)1}=0,\quad N^{(i)}_{(1)j}={1\over 2}\left[\xi^i_{\Vert j}-
\varphi^{ir}\xi^s_{\Vert r}\varphi_{sj}\right].
$$

ii) The torsion {\em\bf T} of the generalized Cartan canonical connection,
induced  by the differential system $(O)$, is determined only by the local
components
$$
R^{(i)}_{(1)jk}=r^i_{jkm}y^m-{1\over 2}\left[\xi^i_{\Vert j\Vert k}-\varphi^
{ir}\xi^s_{\Vert r\Vert  k}\varphi_{sj}\right].
$$

iii)  The electromagnetic components of the metrical time-dependent Lagrange
space of  electrodynamics, induced  by the orbits system $(O)$, have the
expressions
$$
F^{(1)}_{(i)j}={1\over 2}\left[\xi_{i\Vert j}-\xi_{j\Vert  i}\right].
$$
Moreover, they  are governed by the following generalized Maxwell equations:
$$
\sum_{\{i,j,k\}}F^{(1)}_{(i)j\Vert k}=0,\quad
\sum_{\{i,j,k\}}F^{(1)}_{(i)j}\Vert^{(1)}_{(k)}=0.
$$
\end{th}

\subsection{Pfaffian systems}

\hspace{5mm} Let $M=R,\;\varphi_{11}(x)=1$ and $X^{(1)}_{(\alpha)}(t^\gamma,
x)=A_\alpha(t^\gamma)$, where $A_\alpha(t^\gamma)$ is a distinguished 1-form
on $J^1(T,R)$. The $C^2$ solutions of the Pfaffian system
$$
x_\alpha=A_\alpha(t^\gamma),
\leqno{(P)}
$$
are the global minimum  points of the least squares energy functional
$$
{\cal E}_A=\int_T\left[h^{\alpha\beta}(t^\gamma)x_\alpha x_\beta-2A^\alpha(t^
\gamma)x_\alpha+\Vert A\Vert^2(t^\gamma)\right]\sqrt{h}dt,
$$
where $x_\alpha=x^1_\alpha=\partial x^1/\partial t^\alpha,\;A^\alpha
(t^\gamma)=h^{\alpha\mu}(t^\gamma)A_\mu(t^\gamma)$ and $\Vert A\Vert^2(t^\gamma)
=A^\mu(t^\gamma)A_\mu(t^\gamma)$.

\begin{th}
i) The nonlinear  connection on $J^1(T,R)$, induced by the Pfaffian system
$(P)$, has the components
$$
M^{(i)}_{(\alpha)\beta}=-H^\mu_{\alpha\beta}x_\mu,\quad N^{(1)}_{(\alpha)1}
=0.
$$

ii) The only local torsion d-tensor, induced by $(P)$, which does not vanish
is
$$
R^{(1)}_{(\alpha)\beta\gamma}=-H^\mu_{\alpha\beta\gamma}x_\mu.
$$

iii)  The electromagnetic components $F^{(\alpha)}_{(1)1}$, induced by the
Pfaffian system $(P)$, vanish.
\end{th}

\subsection{Continuous groups of transformations}

\hspace{5mm} The  fundamental PDEs system of a transformations group having
the infinitesimal generators $\{\xi_a(x^k)\}_{a=\overline{1,c}}$,
as d-vector fields  on $J^1(T,M)$, is
$$
x^i_\alpha=\sum_{a=1}^c\xi^i_a(x^k(t^\gamma))A_\alpha^a(t^\gamma),
\leqno{(G)}
$$
where $\{A^a(t^\gamma)\}_{a=\overline{1,c}}$ are d-forms on $J^1(T,M)$.

\begin{th}
The $C^2$ solutions of the PDEs system $(G)$ are the global minimum
points of the least squares energy functional
$$
{\cal E}_{\xi A}=\int_T\left[h^{\alpha\beta}\varphi_{ij}x^i_\alpha x^j_\beta-
2A^{\alpha a}\xi_{aj}x^j_\alpha+<A^a,A^b><\xi_a,\xi_b>\right]\sqrt{h}dt,
$$
where
$$
\begin{array}{l}\medskip
A^{\alpha a}(t^\gamma)=h^{\alpha\mu}(t^\gamma)A^a_\mu(t^\gamma),\\\medskip
\xi_{aj}(x^k)=\varphi_{jm}(x^k)\xi^m_a(x^k),\\\medskip
<A^a,A^b>(t^\gamma)=h^{\alpha\beta}(t^\gamma)A^a_\alpha(t^\gamma)A^b_\beta
(t^\gamma),\\
<\xi_a,\xi_b>(x^k)=\varphi_{ij}(x^k)\xi^i_a(x^k)\xi^j_b(x^k).
\end{array}
$$
\end{th}

\begin{th}
i) The canonical nonlinear  connection on $J^1(T,M)$, induced by the PDEs
system $(G)$, is determined by the components
$$M^{(i)}_{(\alpha)\beta}=-H^\mu_{\alpha\beta}x^i_\mu,\quad
N^{(i)}_{(\alpha)j}=\gamma^i_{jm}x^m_\alpha-{A^a_\alpha\over 2}\left[\xi^i_{a
\Vert j}-\varphi^{ir}\xi^s_{a\Vert r}\varphi_{sj}\right].
$$

ii) The local d-torsions which  characterize the geometry of $J^1(T,M)$, attached
to the PDEs system $(G)$, have the expressions
$$
\begin{array}{l}\medskip
\displaystyle{
R^{(i)}_{(\alpha)\beta\gamma}=-H^\mu_{\alpha\beta\gamma}x^i_\mu,}
\\\medskip
\displaystyle{R^{(i)}_{(\alpha)\beta j}={A^a_{\alpha//\beta}\over 2}\left[
\xi^i_{a\Vert j}-\varphi^{ir}\xi^s_{a\Vert r}\varphi_{sj}\right],}\\
\displaystyle{
R^{(i)}_{(\alpha)jk}=r^i_{jkm}x^m_\alpha-{A^a_\alpha\over 2}\left[\xi^i_{a
\Vert j\Vert k}-\varphi^{ir}\xi^s_{a\Vert r\Vert k}\varphi_{sj}\right].}
\end{array}
$$
\end{th}

\begin{th}
i) The components of the electromagnetic field $F$ of the metrical multi-time
Lagrange space attached to the PDEs system $(G)$, are
$$
F^{(\alpha)}_{(i)j}={A^{\alpha a}\over 2}\left[\xi_{ai\Vert j}-\xi_{aj\Vert
i}\right].
$$

ii) The following generalized Maxwell equations, derived  from the PDEs
system $(G)$, hold good:
$$
\sum_{\{i,j,k\}}F^{(\alpha)}_{(i)j\Vert k}=0,\quad
\sum_{\{i,j,k\}}F^{(\alpha)}_{(i)j}\Vert^{(\gamma)}_{(k)}=0.
$$
\end{th}

\subsection{Constrained potentials of physical fields}

\hspace{5mm} Let $(E,\pi,T)$ be a vector bundle of rank $q$ over the Riemannian
temporal manifold $(T,h)$, whose structural group is the compact Lie orthogonal
subgroup $G\subset O_q(R)$. Let $\{\sigma_a\}_{a=\overline{1,q}}$ be
local basis of cross-sections on $E$, and let us consider $g=(g_{ab})_{a,b=
\overline{1,q}}\in G$, the coordinate transformations of  the bundle  $E$.

Let $F\in\Gamma(\Lambda^2\T^{\;*}T\otimes End(E))$ be a fixed $End(E)$-valued 2-form
on $T$, locally expressed by
$$
F=\sum_{1\leq\alpha<\beta\leq p}F_{\alpha\beta}dt^\alpha\wedge dt^\beta,
$$
where the components $F_{\alpha\beta}=(F^a_{b\alpha\beta})_{a,b=\overline{1,q}}$
are skew-symmetric matrices (i. e. , $F_{\alpha\beta}\in L(G)\subset L(O_q(R))=
o(q),\;\forall\;\alpha,\beta\in\{1,\ldots,p\}$). Physically, the cross-section
$F$ is a geometrical model for a given {\it physical field}, in Yang-Mills sense.
For example, $F$ can be one of the gravitational, electromagnetic, strong or
weak nuclear fields \cite{2}.

Let $\nabla\in\Gamma(\Lambda^1\T^{\;*}T\otimes End(E))$ be a $G$-covariant derivative
on $E$, whose locally described by
$$
\nabla=\nabla_\alpha dt^\alpha,
$$
where $\nabla_\alpha=(\nabla^a_{b\alpha})_{a,b=\overline{1,q}}\in L(G)\subset
o(q),\;\forall\;\alpha\in\{1,\ldots,p\}$. From a physical point of view,
the $G$-connection $\nabla$ can be viewed like a possible potential of the
physical field $F$. We recall that a {\it potential of the physical
field $F$}, in the sense of Yang-Mills, is a $G$-connection $\nabla$ on $E$,
which verifies the local identity
\begin{equation}\label{pot}
F_{\alpha\beta}={\partial\nabla_\beta\over\partial t^\alpha}-
{\partial\nabla_\alpha\over\partial t^\beta}+\nabla_\alpha\nabla_\beta-
\nabla_\beta\nabla_\alpha,\quad\forall\;\alpha<\beta\in\{1,\ldots,
p\}.
\end{equation}

In the sequel, we will show that the $G$-potentials  $\nabla$ of the given
physical field $F$ are harmonic maps on the jet space $J^1(T,M)$, where
$M=\Lambda^1\T^{\;*}T\otimes End(E)$.

To reach this aim, let us  suppose that the set of partial derivatives
$\partial\nabla_\alpha/\partial t^\beta,\;\;\forall\;\alpha\leq\beta$ verifies
the {\it constraints} on $J^1(T,M)$,
$$
{\partial\nabla_\alpha\over\partial t^\beta}=f_{\alpha\beta}(t^\gamma,\nabla_
\mu),\quad\forall\;\alpha\leq\beta,\leqno{(C)}
$$
where $f_{\alpha\beta}(t^\gamma,\nabla_\mu),\;\;\forall\;\alpha,\beta\in\{1,2,
\ldots,p\}$,
are {\it "a priori"} given local components of a  geometrical object  on
$J^1(T,M)$, which, under a multi-temporal reparametrization
$(t^\mu)\leftrightarrow(t^{\mu^\prime})$, transform like the local components
$\partial\nabla_\alpha/\partial t^\beta,\;\;\forall\;\alpha,\beta=
\overline{1,p}$. In this context, from PDE equations $(3.1)$ we obtain
$$
{\partial\nabla_\alpha\over\partial t^\beta}=\nabla_\alpha\nabla_\beta-
\nabla_\beta\nabla_\alpha+f_{\beta\alpha}(t^\gamma,\nabla_\mu)-
F_{\alpha\beta}(t^\gamma),\quad\forall\;\alpha>\beta.
$$
Consequently, the PDE equations $(3.1)$, constrained by $(C)$, can be
rewritten in the  equivalent  PDEs system form on $J^1(T,M)$,
$$
{\partial\nabla_\alpha\over\partial t^\beta}={\cal F}_{(\alpha)(\beta)}(t^\gamma,
\nabla_\mu),\quad\forall\;\alpha,\beta=\overline{1,p},\leqno{(CPS)}
$$
where
$$
{\cal F}_{(\alpha)(\beta)}(t^\gamma,\nabla_\mu)=
\left\{\begin{array}{ll}\medskip
f_{\alpha\beta}(t^\gamma,\nabla_\mu),&\forall\;\alpha\leq\beta\\
\nabla_\alpha\nabla_\beta-
\nabla_\beta\nabla_\alpha+f_{\beta\alpha}(t^\gamma,\nabla_\mu)-
F_{\alpha\beta}(t^\gamma),&\forall\;\alpha>\beta.
\end{array}\right.
$$
Obviously, the PDEs system $(CPS)$ is equivalent to
$$
{\partial\nabla^a_{b\beta}\over\partial t^\alpha}={\cal F}^a_{b(\alpha)(\beta)}
(t^\gamma,\nabla_\mu),\quad\forall\;\alpha,\beta=\overline{1,p},\quad\forall
\;a,b=\overline{1,q}.\leqno{(CPS^\prime)}
$$

Taking into account the Riemannian structure induced from $T$ on $\Lambda^1
\T^{\;}*T$,
$$
<dt^\alpha\wedge dt^\beta,dt^\mu\wedge dt^\nu>=h^{\alpha\mu}h^{\beta\nu}-
h^{\beta\mu}h^{\alpha\nu},
$$
and
the well known Riemannian metric on the Lie  algebra $o(q)$,
$$
<A,B>=-Tr(AB^t), \;\;\forall\;A,B\in o(q),
$$
we can endow the vector bundle $M=\Lambda^1\T^{\;*}T\otimes End(E)\to T$ with a
natural Riemannian structure, on each fibre.

Now, applying the ideas of this paper to the PDEs system $(CPS^\prime)$, we
obtain the following  important result:

\begin{th}
Every $G$-potential $\nabla$ of the  physical field $F$, which verifies the
constraints $(C)$, is a harmonic map on $J^1(T,M)$,
being global minimum point of the least squares Lagrangian ${\cal L}=L
\sqrt{h}$, where the Lagrangian function $L$ is defined  by
$$
L=h^{\alpha\beta}h^{\mu\nu}\left\{{\partial\nabla^b_{a\alpha}\over\partial
t^\mu}-{\cal F}^b_{a(\alpha)(\mu)}\right\}\left\{{\partial\nabla^a_{b\beta}
\over\partial t^\nu}-{\cal F}^a_{b(\beta)(\nu)}\right\}.
$$
\end{th}
\addtocounter{rem}{1}
{\bf Remarks \therem} i) The least squares Lagrangian function $L$ from
the previous theorem  is more general than the Yang-Mills Lagrangian
function, whose expression is
$$
YM=h^{\alpha\beta}h^{\mu\nu}F^b_{a\alpha\beta}F^a_{b\mu\nu}.
$$
Consequently, the Euler-Lagrange equations of the least squares Lagrangian
${\cal L}$ may be called  the {\it jet-extension of the Yang-Mills equations}.
Moreover, the $G$-connections, that verify the PDEs system  $(CPS^\prime)$,
may be called {\it jet-harmonic Yang-Mills connections}.

ii) The use of the least squares Lagrangian in the study of PDEs systems,
allows  a natural approach of the multi-time geometric dynamics, in
terms of computer simulation \cite{13},  \cite{19}. Via the multi-time
geometric dynamics and its
computer simulation induced by a given PDEs system  of order one, it is
possible to obtain  an original geometrical classification of these systems.
This study of classification is in our research group attention.
\medskip\\
{\bf Open problem.}  What is the real physical meaning of our geometrical
theory?\medskip\\
{\bf Acknowledgments.} The authors express their sincere thanks to Professors
V. Balan, V. Ob\u adeanu and D. Opri\c s for their critical reading of an earlier
version of this paper.

\begin{center}
University POLITEHNICA of Bucharest\\
Department of Mathematics I\\
Splaiul Independentei 313\\
77206 Bucharest, Romania\\
e-mail:  udriste@mathem.pub.ro\\
e-mail: mircea@mathem.pub.ro\\
\end{center}

\end{document}